\documentclass[aps,prl,twocolumn,groupedaddress,showpacs]{revtex4-1}

\usepackage{graphicx}
\usepackage{hyperref}
\usepackage{amsmath,amssymb}

\newcommand{\xx}{\ensuremath{\mathbf{x}}}
\newcommand{\kk}{\ensuremath{\mathbf{k}}}

\newcommand{\superscript}[1]{\ensuremath{^{\text{#1}}}}

\newcommand{\nd}[0]{\superscript{nd}}

\begin{document}


\title{Phyllotaxis, Pushed Pattern-Forming Fronts and Optimal Packing}


\author{Matthew Pennybacker}
\author{Alan C. Newell}
\affiliation{University of Arizona, Department of Mathematics, Tucson, Arizona 85721, USA}


\date{\today}

\begin{abstract}
We demonstrate that the pattern forming partial differential equation derived from the auxin distribution model proposed by Meyerowitz, Traas and others gives rise to all spiral phyllotaxis properties observed on plants. We show how the advancing pushed pattern front chooses spiral families enumerated by Fibonacci sequences with all attendant self similar properties and connect the results with the optimal packing based algorithms previously used to explain phyllotaxis. Our results allow us to make experimentally testable
predictions.
\end{abstract}


\maketitle

\section{Introduction}

Using a model derived from the pioneering ideas of  Meyerowitz, Traas et al \cite{Jonsson2006,*Reinhardt2003} (MT) on the role of PIN 1 proteins in creating an instability of uniform auxin concentrations near a plant's shoot apical meristem (SAM), this Letter reports on stunning new results which lend credence to the view that almost all of the features of phyllotactic configurations are the result of a pushed pattern forming front whose origin is the MT instability. The front leaves in its wake either whorls or Fibonacci spirals.
The patterns we observe exhibit all known self-similar properties, reveal some new invariants and reproduce the well known van Iterson diagrams associated with discrete algorithms (Levitov, Douady and Couder, Atela, Gol\'{e} and Hotton \cite{Levitov1991,*Douady1992,*Atela2003}) which reflect optimal packing strategies based on ideas and observations of Hofmeister and Snow and Snow \cite{Hofmeister1868,*Snow1952}. Further, the location of the maxima of the auxin fields coincide very closely with the point configurations generated by the discrete algorithms, which suggests that pattern forming systems may be a new tool for addressing optimal packing challenges. In short, instability generated patterns are the mechanism by which plants and other organisms can pursue optimal strategies.

The equation we use \cite{Newell2008,*Newell2008b} for the auxin concentration fluctuation about its mean derives from the continuum approximation to a discrete cell by cell description proposed in \cite{Jonsson2006}.  The resulting PDE is strikingly similar to those found in many pattern forming systems which suggests that Fibonacci spiral patterns should be observable in many physical contexts. The instability of the uniform auxin concentration state which gives rise to the pattern is primarily due to  reverse diffusion. In most equilibrium situations, inhomogeneities in chemical concentrations are smoothed by ordinary diffusion. But in plants the situation is not an equilibrium one. PIN 1 proteins in the interiors of cells move under the influence of an auxin gradient to the cell walls where they orient so as to drive auxin in the direction of its gradient. When the effect is sufficiently strong, the net diffusion is negative and an instability with a preferred length scale occurs.  As the shapes initiated by the instability grow, nonlinear interactions determine the preferred planforms. The resulting PDE for the auxin fluctuation concentration $u(\xx,t)$ turns out to be very close to a gradient flow, and is given by
\begin{figure}
\includegraphics[width=\columnwidth]{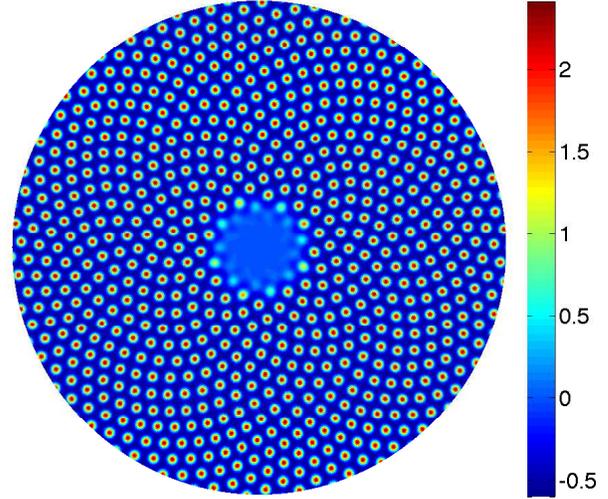}
\caption{A pseudocolor plot of $u(\mathbf{x},t)$ on the inner section $r < 89$ of a pattern initiated at $r = 233$ with parastichy numbers $(89,144)$. A movie may be found at \url{http://math.arizona.edu/~pennybacker/media/sunflower/}.\label{fig:sunflower}}
\end{figure}
\begin{multline}\label{eq:sh}
\frac{\partial u}{\partial t} = -\frac{\delta \mathcal{E}}{\delta u}\\
  = \mu u - \left(\nabla^2 +1 \right)^{\!2}\!u - \frac{\beta}{3}\!\left(|\nabla u|^2 + 2u\nabla^2 u\right) - u^3
\end{multline}
\begin{equation}\label{eq:energy}
\mathcal{E}[u] = -\int \frac{\mu}{2}u^2-\frac12\! \left(u+\nabla^2u\right)^{\!2}-\frac{\beta}{3}u|\nabla u|^2-\frac14 u^4
\end{equation}
where $\xx = (x, y)$ or $(r, \theta)$ are horizontal coordinates on the tunica (plant skin) surface, $t$ is time, the most linearly unstable wavelength is $2\pi$, and $\mu$ and $\beta$ are the coefficients of the linear growth and quadratic terms respectively.

\section{Methods and Data}

We study solutions of \eqref{eq:sh} in many geometries defined by surfaces of revolution $r = r(z)$, $r'(z) = 0$ on a cylinder and $r'(z) = 1$ on a disc, but in this Letter we focus on what we call the sunflower situation in a disc geometry. Sunflowers are formed in two stages. In the first, flowers are initiated in a generative annulus surrounding the plant's SAM and, as the plant grows, the generative annulus and the region of phylla move further out, and their configurations evolve into spiral families. At a certain point, however, the center region consisting of mushy undifferentiated cells begins to undergo a phase transition and, from the outside in, create florets or seeds. In this process, the seeds are laid down annulus by annulus by an inwardly moving front and the particular pattern which forms at a given radius stays at that radius. This means that any optimal packing property which the pattern manifests when it is first laid down remains visible.The diameter of the plant when the process is completed (over a time scale of several days) is of the order of millimeters. Thereafter the plant grows adiabatically until it reaches its mature diameter of 10-15cms. To simulate this situation, we initiate a spiral pattern with parastichies $(m,n)$  in a circle of radius $m+n$ for values of $(m,n)$ ranging from $(55,89)$ to $(89,144)$ and allow the pattern to propagate inwards according to \eqref{eq:sh} (See Fig. \ref{fig:sunflower}). We then analyze the resulting field $u(\xx,t)$.

The numerical algorithm used is nontrivial and is informed by the pioneering recent works of Furihata and Matsuo \cite{Furihata2010}. The key idea is to use the gradient property of \eqref{eq:sh}, by rewriting \eqref{eq:energy} directly in discrete form, and choosing the expression for the discrete gradient which preserves the main dissipation property $d\mathcal{E}/dt < 0$ and ensures both stability and accuracy. For most choices of a discrete approximation for the continuum expression for $\delta \mathcal{E}/\delta u$, this essential property would be compromised. Details are given in \cite{Newell2012,*Pennybacker2013}.

\begin{figure}
\includegraphics[width=\columnwidth]{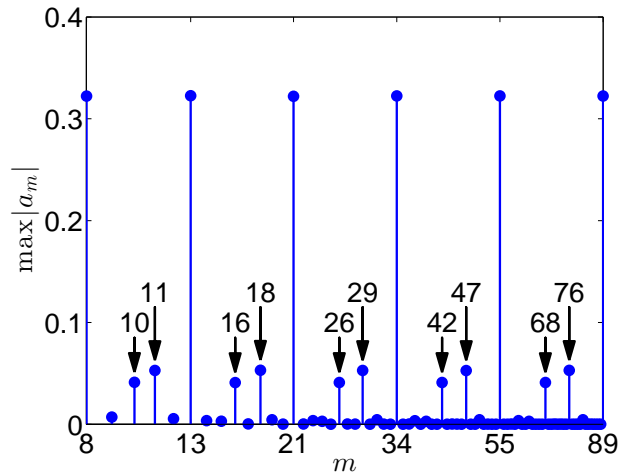}
\caption{The maximum values of the amplitude of all integer modes $8 \leq m \leq 89$. Note that all Fibonacci amplitude maxima are equal, and their 2\nd\ harmonics and ``irregular'' sequence amplitudes are much smaller.\label{fig:amplitude_max}}
\end{figure}

\begin{figure}
\includegraphics[width=\columnwidth]{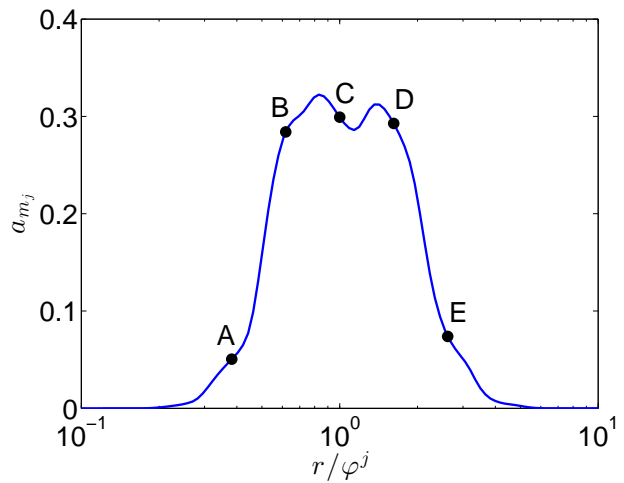}
\caption{The invariant amplitude curve for amplitudes $a_{m_j}$ with $\{m_j\} = \{1,2,3,5,8,13,21,\dots\}$. The horizontal scaling emphasizes the self-similar property $a_{m_j}(r) = a_{m_{j+1}}(\varphi r)$. A movie may be found at \url{http://math.arizona.edu/~pennybacker/media/amplitude/}.\label{fig:amplitude}}
\end{figure}

Guided by previous analytical results at near-onset conditions, we understood several key points.  First is the fact that Fibonacci spirals, and indeed whorl structures, are very much a consequence of the presence of a sign reversal asymmetry and the fact that the pattern is laid down
by a moving front annulus by annulus in an environment of slowly
changing metric. The asymmetry gives hexagonal lattices a special role in planar patterns because modes $\exp(i\kk_m\cdot x)$, having $\kk_m = k_0$, with wavevectors 120 degrees apart reinforce
each other via quadratic interactions. Likewise, in a circular
geometry, for patterns laid down annulus by annulus, we observe triads of modes having wavevectors $\kk_m = (m,\ell_m)$, $\kk_n = (n,\ell_n)$, $\kk_{m+n} = \kk_m+\kk_n$ with shapes $exp(-i\int \ell_m dr-im\theta)$, a good
approximation for $r$ large, $m/r$ finite. The signs of the radial wavenumbers
alternate. These modes can reinforce each other via quadratic
interactions at certain radii where they have almost equal
amplitudes and energetically preferred wavenumbers, with lengths given by $\kk_m^2 = \ell_m^2+m^2/r^2$. But, as the radius
increases, the lengths of the wavevectors begin to deviate from their optimal values.
Moreover, the energetically preferred radial wavenumbers $(\ell_m(r), \ell_n(r))$
 evolve along a computable locus. As a consequence, the first
three of the new quadratically generated modes, $\kk_{2m}$, $\kk_{2n}$,
$\kk_{2m+n}$ move away from, whereas the last wavevector, $\kk_{m+2n}$,
moves towards the critical circle. This results in the amplitudes of
the first three modes being slaved to, and much less than, the amplitude
of the last (see Fig. \ref{fig:amplitude_max}). Thus, in situations where the pattern is laid
down annulus by annulus, the quadratic interactions select from the
set of all integers only that subset obeying Fibonacci rules.

These considerations suggest that the quantities we should monitor are the amplitudes and radial wavenumbers of the signal $u(\xx,t) = \sum_j u_{m_j}(r,t)\exp(-im_j\theta) + (*)$, $\{m_j\} \subset \mathbb{Z}$, which we obtain by writing $u_{m_j}(r,t)$ as $a_{m_j}(r,t) \exp(i \phi_{m_j})$ . At any given $r$ and $t$, the amplitudes are given by the sequence $\{a_{m_j}\}$ and the radial wavenumbers $\ell_{m_j}$ by $d \phi_{m_j}/dr$. We also monitor the front speed $\nu$, which is not constant but varies in a periodic fashion over distances separated by $\varphi$, the golden number. Because the chosen pattern combines several modes, one has to ask if they can all propagate as a synchronous front. There are two kinds of fronts \cite{VanSaarloos2003}. The first kind are pulled fronts whose properties are determined by conditions in the unstable state ahead of the front. The second category, to which Fibonacci spirals belong, are pushed fronts whose speeds and steepnesses exceed those of the pulled front and whose characteristics are determined by conditions behind the front. This is essential for Fibonacci patterns because it ensures both that the modes participating in the pattern structure all propagate in synchrony and introduces the important ingredient of bias into pattern choice. Namely the choice of pattern emerging in the present generative annulus is influenced by the pattern already laid down in the previous one. This is important because the free energy landscape has many minima. The minimum which is chosen is the one for which the previous state of the system lies in its basin of attraction. The minima corresponding to Fibonacci spirals are metastable states, continuous in $r$. They are not the deepest minima in the landscape. Whereas the metastable Fibonacci  pattern state  is long lived, it will eventually tunnel its way into a lower energy minimum if one is available on the landscape. And eventually of course, the pattern will occupy a region large enough so that it will think it is on a plane. Very slowly, it will develop defects which will lead to patches of pure hexagonal planforms mediated by defects. But coarsening takes long time. Moreover, in plants there may be a mechanism which roots the angular position of a primordium once it begins to grow.

\section{Results}

\begin{figure}
\includegraphics[width=\columnwidth]{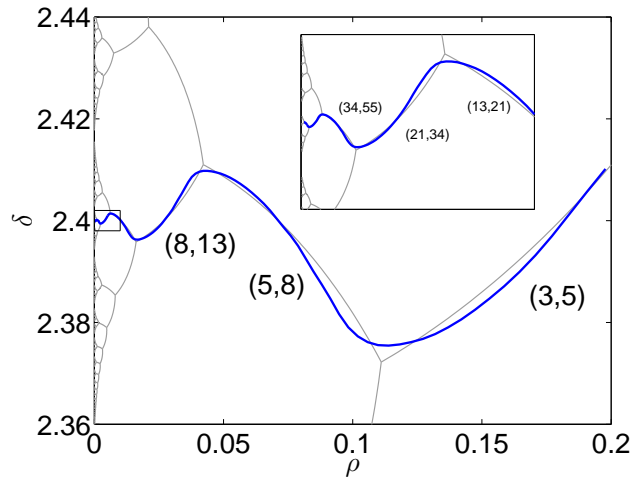}
\caption{The rise $\rho$ and divergence angle $\delta$ given by the local lattice structure at each radius. The shaded lines are the van Iterson tree, with selected parastichy pairs indicated. Inset is detail of the data for small $\rho$. A movie may be found at \url{http://math.arizona.edu/~pennybacker/media/divergence/}.\label{fig:divergence}}
\end{figure}

\begin{figure}
\includegraphics[width=\columnwidth]{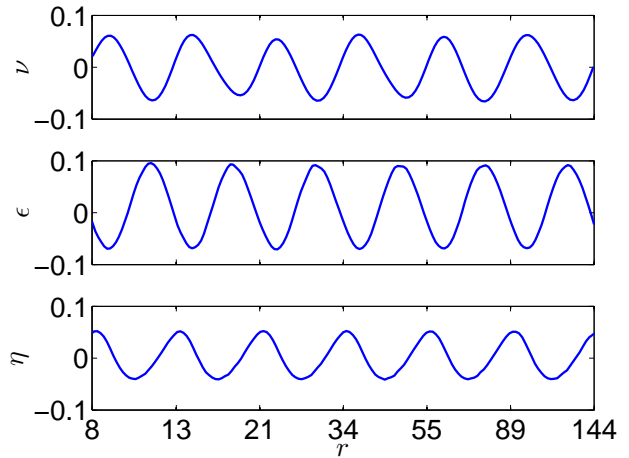}
\caption{The front speed $\nu$, the local energy $\epsilon$, and the local packing efficiency $\eta$. The vertical axis has been rescaled to indicate relative variation from the mean value of each of these quantities.\label{fig:energy}}
\end{figure}

\paragraph{Simulation data.}
In Fig. \ref{fig:sunflower}, we show the simulated sunflower head. The transitions between pairs of parastichy numbers are smooth. For a wide range of $\mu$, $\beta$, the result is the same; we use $\mu = 1\times 10^{-3}$ and $\beta = 3$ for all of the data that follow.  The pattern makes its transitions from one almost hexagonal state to another via intermediate two-mode structures, whose fields look rhombic. In Fig. \ref{fig:amplitude_max}, we plot the maximum amplitudes (over the sunflower head) of all integer modes $8 \leq m \leq 89$. The Fibonacci modes $\{m_j\} = \{8,13,\dots,89\}$ are dominant, although the 2\nd\ harmonics $\{2m_j\}$ and irregular Fibonacci modes $\{m_{j-2}+m_j\}$, passively slaved, have significant albeit small amplitudes. They measure the small difference between a ``perfect'' and a pattern formed Fibonacci lattice. The amplitudes of all other integer modes are insignificant. The choice of the ``winning'' subsequence of dominant modes only depends on the initial condition. Different Fibonacci-like sequences of dominant wavenumbers even inherit the same self-similar properties discussed in the following remarks. We also check in every case the self-similar property that $a_{m_{j+1}}(r\varphi) = a_{m_j}(r)$, and that the radial wavenumbers satisfy a similar relation $\ell_{m_{j+1}}(r\varphi) = -\ell_{m_j}(r)$, where $\varphi = (1+\sqrt{5})/2$ is the golden number.

\paragraph{The field amplitudes.}
In Fig. \ref{fig:amplitude}, we show the new ``amplitude invariant'' which has no parallel in the discrete algorithms which only generate point configurations. We plot, as functions of $r/\varphi^j$ on a logarithmic scale, the locus of all amplitudes $\{a_{m_j}\}$, $\{m_j\} = \{3,5,8,\dots,144\}$, as $r$ decreases from $233$ to $8$. At $r \approx 144$, the configuration is almost hexagonal with $a_{144}$, $a_{89}$, $a_{55}$ occupying the positions $B$, $C$, $D$ respectively. As $r$ decreases, the amplitude $a_{144}$ will very quickly decrease leaving $a_{55}$ and $a_{34}$ dominant (and almost unchanged) for a short interval until the amplitude $a_{21}$ rises quickly to form the next hexagonal configuration at $r \approx 89$. As $r$ continues to decrease, the pattern repeats with $a_{89}$, $a_{55}$, $a_{34}$ and $a_{21}$ replacing $a_{144}$, $a_{89}$, $a_{55}$ and $a_{34}$. The graph of the locus of $\{a_{m_j}\}$ in three dimensions $(a_{m_{j+1}}, a_{m_j}, a_{m_{j-1}})$ is given in \cite{Newell2012} and shows a common homoclinic orbit joining the origin to itself with six legs on which there are two-mode dominated transitions from hexagons from circumferential wavenumbers $m_{j+1}$, $m_j$, $m_{j-1}$ to $m_j$, $m_{j-1}$, $m_{j-2}$.

\paragraph{The lattice configurations.}
The overall arrangement of phylla on the sunflower head is not a fixed spiral lattice, but a slowly varying one. A description of this variation may be found by assigning to each radius a ``local'' lattice. On the circle of radius $r$, we compute the amplitudes and wavevectors of the dominant modes, which give the orientations of the parastichies and the shapes of the phylla. In fact, we can measure from the simulation data the set of wavenumbers $\{\ell_m\}$ for all integers $m$. The curves defined, as $r$ varies, by $(\ell_m,\ell_n)$ for all consecutive pairs, which belong to some Fibonacci subsequence, of dominant circumferential wavenumbers $m,n$ are identical.
In Fig. \ref{fig:divergence}, we choose to present one manifestation of this property, which allows a direct comparison with the van Iterson diagram,  by performing a trivial change of coordinates onto the cylinder of radius $r$, taking the radial coordinate to the axial coordinate. This preserves the differential structure in the neighborhood of the circle. Supposing that the two dominant modes have circumferential wavenumbers that are coprime, the resulting lattice has maxima only at distinct and regularly spaced axial positions. The axial distance $\rho$, normalized by the radius of the cylinder, and the angular difference $\delta$ between two consecutive maxima, often called rise and divergence angle, uniquely identify the local lattice. This data is overlaid on a tree diagram, due originally to van Iterson, that identifies lattices which are rhombic, having at least four nearest neighbors of each point equidistant.  Rhombic lattices are optimally packed, in that, for a given rise or divergence angle, they have the largest fraction of their area covered by identical circles centered at the lattice points. Hexagonal lattices are a special case, having six equidistant neighbors. It is clear that each of our local lattices is nearly rhombic, and moreover, that the limiting value of the divergence angle is the golden angle $2\pi/\varphi^2 \approx 2.4$.

\paragraph{Self-similar quantities.} In Fig. \ref{fig:energy}, we show the variation of three important quantities that are unchanged by the transformation $r \rightarrow r\varphi$. The velocity $\nu$ is the rate at which the front propagates inward. The local energy $\epsilon$ is the energy density of the local lattice, computed with the amplitudes and wavevectors measured at each radius. The local packing efficiency $\eta$ is the largest fraction of the local cylinder, representing the pattern in a narrow annulus at radius $r$, that can be covered by identical circles centered at maxima of the local lattice. When the pattern is nearly hexagonal, the local energy is low and both the front velocity and local packing efficiency are high. When the pattern is two-mode dominated, but less hexagonal, the opposite is true. Recent experimental work \cite{Hotton2006} on masking developing seed patterns in sunflowers suggests that it may be possible to measure $\nu$ and check that it varies as predicted. 

Whorl structures with wavevectors $(\pm\ell,m)$, $(0,2m)$ such as the
decussate $m=2$ only exist in finite $r$
intervals. In \cite{Newell2012}, we show how a decussate
makes a transition to a spiral structure via an Eckhaus-like
instability. We examine a range of such transitions. Finally,
we remark that our work suggests that Fibonacci patterns are
universal, long lived but not infinitely long lived, in pattern
forming planar systems for which the plane is tiled annulus by
annulus. More generally, it
will be interesting to explore the nature of patterns laid down by a
front near regions in which the metric changes slowly.

\begin{acknowledgments}
This work was supported by NSF grant DMS 0202440.
\end{acknowledgments}

\bibliography{references}

\end{document}